\documentclass[11pt]{article}

\usepackage[dvips]{color}
\usepackage[dvips]{graphicx}

\usepackage{charter}
\usepackage[T1]{fontenc}
\usepackage{textcomp}

\usepackage{calrsfs}

\newcommand{\runninghead}[1]{\gdef\RH{#1}}\newcommand{\RH}{}
\newcommand{\msc}[1]{\gdef\MSC{#1}}\newcommand{\MSC}{}
\renewcommand{\title}[1]{\gdef\TT{#1}}\newcommand{\TT}{}
\renewcommand{\author}[1]{\gdef\AU{#1}}\newcommand{\AU}{}
\newcommand{\address}[1]{\gdef\ADR{#1}}\newcommand{\ADR}{}
\renewcommand{\date}[1]{\gdef\DD{#1}}\newcommand{\DD}{}
\newcommand{\version}[1]{\gdef\VER{#1}}\newcommand{\VER}{}

\renewcommand{\maketitle}{
\par \noindent {\footnotesize \upshape Running head: \RH \hfill \DD \par
\noindent Math.\ Subj.\ Class. (2000): \MSC \hfill\textbf{\VER}} 
\vspace{3cm}
\begin{center} {\LARGE \normalfont \TT} \par \medskip \AU \end{center} \vspace{1cm}
}

\newcommand{\finalinfo}{
\bigskip \noindent {\small \upshape \ADR} }

\usepackage{sectsty}

\allsectionsfont{\normalsize \bfseries \nohang \centering}
\makeatletter
\def\@seccntformat#1{\csname the#1\endcsname. \ }
\makeatother

\usepackage[thmmarks]{ntheorem}

\theoremnumbering{Alph}

\theoremstyle{plain}
\theoremheaderfont{\normalfont\bfseries}
\theorembodyfont{\normalfont\itshape}
\theoremseparator{.}

\newtheorem{introthm}{Theorem}

\theoremstyle{nonumberplain}
\theoremheaderfont{\normalfont\bfseries}
\theorembodyfont{\normalfont\upshape}
\theoremseparator{.}

\theoremstyle{nonumberplain}
\theoremheaderfont{\small\bfseries}
\theorembodyfont{\small\upshape}
\theoremseparator{.}
\newtheorem{ack}{Acknowledgement}

\theoremnumbering{arabic}

\theoremstyle{change}
\theoremheaderfont{\normalfont\bfseries}
\theorembodyfont{\normalfont\itshape}
\theoremseparator{.}

\newtheorem{thm}{Theorem}[section]
\newtheorem{prop}[thm]{Proposition}

\newtheorem{lem}[thm]{Lemma}

\theoremstyle{change}
\theoremheaderfont{\normalfont\bfseries}
\theorembodyfont{\normalfont\upshape}
\theoremseparator{}

\newtheorem{se}[thm]{}

\theoremstyle{change}
\theoremheaderfont{\normalfont\bfseries}
\theorembodyfont{\normalfont\upshape}
\theoremseparator{.}

\newtheorem{ex}[thm]{Example}
\newtheorem{rem}[thm]{Remark}

\theoremstyle{nonumberplain}
\theoremheaderfont{\normalfont\itshape}
\theorembodyfont{\normalfont\upshape}
\theoremseparator{.}
\newcommand{\eop}{\hbox{\rlap{$\sqcap$}$\sqcup$}}
\theoremsymbol{\eop}
\newtheorem{proof}{Proof}

\usepackage{amssymb}

\def \smallmat #1 #2 #3 #4 { \protect \scriptsize \left( \hspace{-0,15cm}
                              \begin{array}{cc} {#1} & \hspace{-0,15cm} #2 \\ #3 & \hspace{-0,15cm} #4 \end{array}
                              \hspace{-0,15cm} \right) }

\def \tinymat #1 #2 #3 #4 { \protect \tiny \left( \hspace{-0,2cm}
                              \begin{array}{cc} {#1} & \hspace{-0,2cm} #2 \\ #3 & \hspace{-0,2cm} #4 \end{array}
                              \hspace{-0,2cm} \right) }                                     
\def \mat #1 #2 #3 #4 { \left( \begin{array}{cc} {#1} & #2 \\ #3 & #4 \end{array} \right) }

\def \vec2 #1 #2  { \protect \scriptsize \left( \hspace{-0,15cm}
                              \begin{array}{c} {#1} \\ {#2}  \end{array}
                              \hspace{-0,15cm} \right) }

\usepackage{bbm}
\def\1{\mathbbm{1}}
\def\IM{\mathop{\mathrm{im}}}
\def\diag{\mathop{\mathrm{diag}}}

\def\Z{\mathbf{Z}}
\def\Q{\mathbf{Q}}
\def\R{\mathbf{R}}

\def\P{\mathbf{P}}

\def\ker{\mathop{\mathrm{ker}}\nolimits}

\def\GL{\mathop{\mathrm{GL}}\nolimits}
\def\PSL{\mathop{\mathrm{PGL}}\nolimits}
\def\an{\mathrm{an}}

\def\T{{\mathcal T}}
\def\cT{{\mathcal T}}

\def\cO{{\mathcal O}}
\def\ord{\mathop{\mathrm{ord}}}

\usepackage[OT2,T1]{fontenc}
\DeclareFontFamily{U}{wncy}{}
\DeclareFontShape{U}{wncy}{m}{n}{<->wncyr10}{}
\DeclareSymbolFont{mcy}{U}{wncy}{m}{n}
\DeclareMathSymbol{\Ch}{\mathord}{mcy}{"51}

\usepackage[british]{babel}

\begin{document}

\runninghead{Graph $C^*$-algebras}
\version{version 2.0}
\date{July, 6th, 2006} 
\msc{ 05C50, 20E08, 46L80}
\title{On the $K$-theory of graph $C^*$-algebras}
\author{\textit{by} Gunther Cornelissen, Oliver Lorscheid \textit{and} Matilde Marcolli} 
 
\address{(gc, ol) Mathematisch Instituut, Universiteit Utrecht, Postbus 80.010, 3508 TA Utrecht, Nederland, \\ email: {\tt cornelis@math.uu.nl, lorschei@math.uu.nl} 

\medskip

\noindent (mm) Max-Planck-Institut f\"ur Mathematik, Vivatsgasse 7, 53111 Bonn, Deutschland, \\ email: {\tt marcolli@mpim-bonn.mpg.de}} 

\maketitle

\begin{abstract} 
\noindent We classify graph $C^*$-algebras, namely, Cuntz-Krieger algebras associated to the Bass-Hashimoto edge incidence operator of a finite graph. This is done by a purely graph theoretical calculation of the $K$-theory and the position of the unit therein.
\end{abstract}

\section*{Introduction} 

Let $\Ch$ denote a finite (multi-)graph with $m$ (geometrical) edges. Note that $\Ch$ might have loops, multiple edges and sinks.  Consider the (oriented) graph $\Ch^+$, whose vertices equal those of $\Ch$, and whose edges consist of all edges of $\Ch$ with both possible orientations. We denote the two edges of $\Ch^+$ corresponding to an edge $e$ of $\Ch$ by $e$ and $\bar{e}$. Let $o$ and $t$ denote the origin, respectively, terminal vertex of an oriented edge $e$ of $\Ch^+$. On the space of edges of the oriented graph $\Ch^+$, we consider the operator 
$$ T : E \Ch^+ \rightarrow E \Ch^+ \ : \ e \mapsto \sum_{t(e)=o(e')}  e'-\bar{e}. $$
This operator is represented by a $2m \times 2m$-matrix $A$, whose entries are in $\{0,1\}$. The operator $T$ was considered by Hashimoto \cite{Hash} and Bass \cite{Bass} in connection with their study of the Ihara zeta function of a graph (see also: Stark and Terras \cite{ST}). 

Let $\cO_{\Ch}$ denote the Cuntz-Krieger algebra (\cite{CK}) associated to the matrix $A$. This is the universal $C^*$-algebra generated by $2m$ partial isometries $\{S_i\}_{i=1}^{2m}$ with orthogonal range projections, subject to the relations $S_i^* S_i = \sum A_{ij} S_j S_j^*$. 

We want to classify Cuntz-Krieger algebras that arise in such a way, up to the three natural notions of (strong) Morita equivalence (denoted $\sim$), stable isomorphism (also denoted $\sim$ in this paper, since the notions are equivalent in our setting) and strict isomorphism (denoted $\cong$). Although it seems at first that quite a few algebras are possible, we will see that this turns out not to be the case. 

The expert in operator algebras might immediately see the picture: Kum\-jian and Pask (\cite{KP}) have proven that $\cO_\Ch$ is Morita equivalent to a boundary operator algebra. More specifically, let $\cT$ denote the universal covering tree of $\Ch$. Then $\Ch=\cT/\Gamma$ for $\Gamma$ a free group of rank $g=$ the first Betti number of $\Ch$, and $\cO_\Ch \sim C^*(\partial \cT) \rtimes \Gamma$, where the right hand side algebra, which we call a \emph{boundary operator algebra}, only depends on $g$. As a matter of fact, since a result of R{\o}rdam (\cite{Rordam}) implies that $K_0$ is a full invariant for the stable isomorphism class of Cuntz-Krieger algebras and stable isomorphism and Morita equivalence are the same for such algebras, this result can be used to compute the $K$-theory of $\cO_\Ch$ by calculating it for one example graph $\Ch$, and this is done by Robertson in \cite{Rob} in an operator algebraic way. 

Note that R{\o}rdam  has also proven that $K_0(\cO_\Ch)$ together with the position of the unit in $K_0(\cO_\Ch)$ is a full  invariant for the strict isomorphism type of $\cO_\Ch$. 

Our aim in this paper is twofold. First, we want to compute the $K$-theory of $\Ch$ directly in a graph theoretical way, thus avoiding the boundary operator algebra. Our proof will involve relating the Smith Normal Form reduction of $1-A$ to contraction of non-loops in the graph. The proof of our second result will rely heavily on this technique. Although knowing $K_0$ for a Cuntz-Krieger algebra allows one to compute $K_1$, we also exhibit a direct relation between $K_1$ and graph homology. The result is

\begin{introthm} \label{thm1} Let $\Ch$ denote a finite graph with first Betti number $g \geq 1$. Then $$K_0(\cO_\Ch) \cong \Z^g \oplus \Z/(g-1)\Z. $$ Furthermore, $K_1(\cO_\Ch)$ is naturally isomorphic to $Z_1$, the space of cycles on $\Ch$, so $K_1(\cO_\Ch) \cong \Z^g$. 
Hence for $g \geq 2$, the Morita equivalence and stable equivalence type of $\cO_\Ch$ only depends on $g$. 

\end{introthm}

The result can now be used, conversely, to compute the $K$-theory of boundary operator algebras. One can also deduce some results about the operator $T$ from this theorem. In particular, $1-T$ has kernel $K_1(\cO_\Ch)$, $\Z^{|E\Ch^+|}/\IM(1-T) \cong K_0(\cO_\Ch)$, and one recovers the computation of the rank of $1-T$ from \cite{Bass} for both $g \geq 2$ and $g=1$.

Secondly, we want to study the strict isomorphism class of $\cO_\Ch$ for a general graph $\Ch$. This problem doesn't seem to have been dealt with in the literature. By R{\o}rdam's classification results, it is known to depend upon the position of the unit in $K_0$. Again from the boundary operator algebra point of view, Robertson showed in \cite{Rob} that the unit of $C^*(\partial \cT) \rtimes \Gamma$ is of exact order $g-1$. This is analogous to  Connes's calculation in \cite{Connestrans}, Cor.\ 6.7, that the class of the identity $\1$ in $K_0(A)$ has order $g-1$, where $A=C(\P^1(\R)) \rtimes \Gamma$, with $\Gamma$  a torsion-free cocompact lattice in $\PSL(2,\R)$, and $g$ the genus of the Riemann surface uniformized by $\Gamma$. 

Here, we use purely graph theoretical considerations to calculate the exact position of the unit in $K_0(\cO_\Ch)$ for a general graph $\Ch$. Let us denote by $(a,b)$ the greatest common divisor of two integers $a$ and $b$. 

\begin{introthm} \label{thm2} 
Let $\Ch$ denote a graph with first Betti number $g \geq 2$. Then the image of the unit of $\cO_\Ch$ in $K_0(\cO_\Ch) \cong \Z^g \oplus \Z/(g-1)\Z$ has (finite) order $$ \frac{g-1}{(g-1,|V\Ch|)}, $$
where $|V\Ch|$ is the number of vertices of $\Ch$. In particular, the class of the unit in $K_0$ is anihilated by the Euler characteristic $g-1$ of $\Ch$. Furthermore, every possible strict isomorphism type of an operator algebra with $K_0$-group of the given type occurs as graph $C^*$-algebra for a stable graph; strict isomorphism is not a homotopy invariant for graphs; and $\cO_\Ch$ is only strictly isomorphic to a boundary operator algebra of genus $g$ if the number of vertices in $\Ch$ is coprime to the Euler characteristic $g-1$ of $\Ch$. 
\end{introthm} 

We hope that the purely combinatorial considerations leading to these theorems will help in understanding the higher dimensional analogues of these results  (for group actions on buildings) as considered by Robertson in \cite{Rob2}, \cite{Rob}.

Our results arose from trying to answer the following question in algebraic geometry: one can associate a spectral triple to a Mumford curve $X$ (cf.\ \cite{ConsaniMarcolli}), in which the operator algebra is the boundary operator algebra for the Schottky group. In this construction, one can replace that algebra by the graph $C^*$-algebra $\cO_\Ch$ for $\Ch$ the stable reduction of $X$. Does this new algebra capture more geometrical information   about $X$? See Remark \ref{stable} below for some more details.

We only consider finite graphs. It might be interesting to extend the methods to row-finite or locally finite graph $C^*$-algebras. 

\begin{ack} 
We thank Jakub Byszewski for his input in section \ref{k1}. The position of the unit in $K_0(\cO_\Ch)$ was guessed based on some example calculations by Jannis Visser in his SCI 291 Science Laboratory at Utrecht University College.
\end{ack}

\section{Stable isomorphism type}\label{classO}

\begin{se}[Set-up] Let $\Ch$ denote a graph (possibly with loops, multiple edges, and sinks) with $m$ (geometrical) edges and of first Betti number (``cyclomatic number'') $g$, viz., $g$ equals the number of independent loops, which equals the number of edges outside a spanning tree. We will mostly be considering the case where $g\geq 2$, but we will comment upon what happens for $g=1$ at the appropriate place. In this section, we will prove Theorem \ref{thm1} from the Introduction.
\end{se}

\begin{se}[The operator $T$] We make a new (oriented) graph $\Ch^+$, whose vertices equal those of $\Ch$, and whose edges consist of all edges of $\Ch$ with both possible orientations. We denote the two edges of $\Ch^+$ corresponding to an edge $e$ of $\Ch$ by $e$ and $\bar{e}$. Let $o$ and $t$ denote the origin, respectively, terminal vertex of an oriented edge $e$ of $\Ch^+$. On the space of edges of the oriented graph $\Ch^+$, we consider the operator 
$$ T : E \Ch^+ \rightarrow E \Ch^+ \ : \ e \mapsto \sum_{t(e)=o(e')}  e'-\bar{e}. $$
This operator is represented by a $2m \times 2m$-matrix $A$, whose entries are in $\{0,1\}$. The operator $T$ was considered by Hashimoto \cite{Hash} and Bass \cite{Bass} in connection with their study of the Ihara zeta function of a graph (see also: Stark and Terras \cite{ST}). Although not strictly necessary for the main argument of this paper, we will comment upon this relation in \ref{graphzeta} below.
\end{se}

\begin{se}[The Cuntz-Krieger algebra] Let $\cO_{\Ch}$ denote the Cuntz-Krieger algebra associated to the matrix $A$. This is the universal $C^*$-algebra generated by $2m$ partial isometries $\{S_i\}_{i=1}^{2m}$ with orthogonal range projections, subject to the relations $$S_i^* S_i = \sum A_{ij} S_j S_j^*.$$ 
 We want to classify Cuntz-Krieger algebras that arise in such a way. Although it seems at first that quite a few algebras are possible, we will see that this turns out not to be the case. 

First some preliminary observations. Since the reduction graph has cyclomatic number $g \geq 2$, it has a vertex of valency $\geq 2$, so $A$ is not a permutation matrix (actually, the only exception would be a graph consisting of a single vertex and a single loop, which corresponds to a Tate elliptic curve). The matrix $A$ is also irreducible, i.e., all entries of $A^m$ for $m$ sufficiently large are positive. This is clear from the interpretation of $A^m_{x,y}$ as the number of paths of length $m$ between the end point of $x$ and the origin of $y$, since $\Ch$ is a finite and connected graph. Therefore, $\cO_{\Ch}$ is a simple algebra (\cite{CK}, Theorem 2.14). By a result of R{\o}rdam (\cite{Rordam}), two different such algebras for graphs $\Ch$ and $\Ch'$ are stably isomorphic (i.e., isomorphic after tensoring with compact operators) if and only if $K_0(\cO_{\Ch}) \cong K_0(\cO_{\Ch'})$. Furthermore, they are isomorphic if and only if there is a group isomorphism between the $K_0$-groups that maps the class of the unit of $\cO_{\Ch}$ to that of $\cO_{\Ch'}$. The image of the unit in the abstract group $K_0(\cO_\Ch)$ up to abstract group automorphisms is called ``the position of the unit'' in the literature. 
\end{se}

\begin{rem} \label{BGR}
The notions of stable isomorphism and (strong) Morita equivalence are equivalent for algebras of the form $\cO_{\Ch}$, since they are separable, see Brown-Green-Rieffel \cite{stablemorita}.
\end{rem}

\begin{se}[$K$-theory computations] We hence start by computing the $K$-theory of $\cO_{\Ch}$. Theorem 5.3 in \cite{CK} says that $K_0(\cO_{\Ch}) = \Z^n /(1-A^t)\Z^n$, where $n$ is the dimension of the matrix $A$. Equivalently, in terms of the operator $T$, 
$K_0(\cO_{\Ch}) = \Z^{E\Ch^+}/\IM (1-T).$ Note that it is irrelevant whether one works with $A^t$ or $A$ in this formula. We now claim the following:
\end{se}

\begin{prop} \label{oli}
If $\Ch$ is a graph with cyclomatic number $g \geq 1$, then $$K_0(\cO_{\Ch}) \cong \Z^g \oplus \Z/(g-1)\Z;$$ in particular, the stable isomorphism type of $\cO_{\Ch}$ only depends on $g$ if $g \geq 2$.
\end{prop}

\begin{rem}
As follows from \cite{Rordam} 4.3, any abelian group can be the $K$-group of a Cuntz-Krieger algebra associated to the \emph{vertex} adjacency matrix of a graph. In light of this, the above result might come as a surprise.
\end{rem}

\begin{proof} The proof has two parts. We use the following notation: for two square matrices $M,N$ of size $n \times n$, we write $M \sim N$ if there exist matrices $X,Y \in \GL_{n}(\Z)$ such that $XMY=N$. We also write $1_n$ for the indentity matrix of size $n \times n$.

\medskip

\begin{figure}[h]
   \centering
   \input{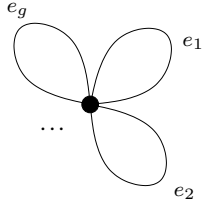} 
   \caption{The genus $g$ graph with one vertex}
   \label{fig1}
\end{figure}

\noindent \textbf{\ref{oli}.1} We start by proving this for the case where $\Ch$ has only one vertex, so $m=g$, see Figure \ref{fig1}. 

The matrix $A-1$ is of the form
$ A-1 = \tinymat B B B B ,$ where $B$ is a $g \times g$ matrix with zeros along the diagonal and $1$ everywhere else.

Obviously, $A-1 \sim \tinymat B 0 0 0 .$ By subtracting the first column from all other in $B$, and then adding all rows to the last, then adding all columns to the first, we find $B \sim \diag(g-1,1,\dots,1)$. It follows that $$\IM(1-A) \cong \IM(\diag(1,\dots,1,g-1,\underbrace{0,\dots,0}_g))=\Z^{2m-g-1} \oplus (g-1) \Z,$$ and the result follows.

\medskip

\noindent \textbf{\ref{oli}.2} We now consider what happens to the operator when we contract an edge.

\begin{quote} \textbf{Claim} \ \emph{If $\Ch'$ is the graph obtained from contracting a single non-loop edge in $\Ch$, and $A'$ and $A$ are the matrices of the respective $T$-operators on these graphs, then $ A-1 \sim \tinymat 1_2 0 0 A'-1 $.} \end{quote}

\noindent Suppose that $\gamma$ is the edge that is contracted. Suppose $\gamma$ and $\bar{\gamma}$ correspond to the first and second row of $A-1$ respectively. Perform the following elementary row operations on $A-1$: add the first row to every row corresponding to an edge $e$ whose terminal point is the origin of $\gamma$; and add the second row to every row corresponding to an edge $e$ whose terminal point is the origin of $\bar{\gamma}$. 

\medskip

\begin{figure}[h]
   \centering
   \input{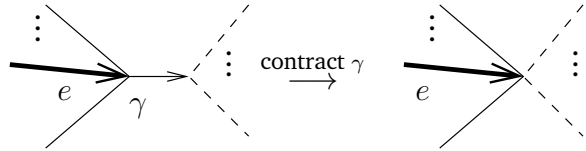} 
   \caption{Contracting the edge $\gamma$}
   \label{fig2}
\end{figure}

Observe the following facts about the resulting matrix: 
\begin{enumerate}
\item Since $\gamma$ is not a loop, the left upper corner is $-1_2$.
\item The cofactor of the first two rows is equal to $A'-1$; because of the transformation that we did, see Figure \ref{fig2}: the outgoing (dashed) edges of $\gamma$ become added to the outgoing edges of $e$.
\item All entries in the first two columns from the third row on are zero; because in the original matrix, there was a $1$ at location $(e,\gamma)$ exactly if $e$ flows into $\gamma$, but then our operation has added the $-1$ from location $(\gamma,\gamma)$ to it --- and similarly for the second row of $(e,\bar{\gamma})$.
\end{enumerate}
Now any possible non-zero entries in the first two rows from the third column on can be removed by adding the first or second column to the corresponding column, and this will not affect any other entry. In the end, we find a matrix as in the Claim above.

\medskip

Given a general graph $\Ch$, we contract all non-loops one after the other. By the claim, we are left with a matrix of the form $1-A \sim \tinymat -1_{2m-g} 0 0 A'-1 $, where $A'$ corresponds to a graph with only loops, i.e. with one vertex. But then \ref{oli}.1 can be applied to $A'-1$, and we find in the end that$$ 1-A \sim \diag(g-1,1,\dots,1,\underbrace{0,\dots,0}_g).$$ As we remarked above, the algebra $\cO_\Ch$ is simple for $g \geq 2$, so $K_0$ classifies it up to stable isomorphism. This finishes the proof of Proposition \ref{oli}. 
\end{proof}

\begin{se}[$K_1$ and graph homology] \label{k1} Since $K_1$ of a Cuntz-Krieger algebra is isomorphic to the torsion free part of $K_0$ (see \cite{Rordam}), the same proof shows that $K_1(\cO_{\Ch}) \cong \ker(1-T) \cong \Z^g$ for $g \geq 2$; and that $K_1(\cO_{\Ch}) \cong \ker(1-T) \cong \Z^2$ for $g=1$.  We now make a digression to interpret this result in terms of graph  homology.

Let $Z_1=H_1(\Ch,\Z)$ denote the space of (integral) cycles on $\Ch$, i.e., the kernel of the boundary map $\Z^{E \Ch} \rightarrow \Z^{V \Ch}.$ We observe that one can decompose the space $\Z^{E \Ch^+}$ as follows:
$$ \Z^{E \Ch^+} = \Z^{S_\Ch} \oplus Z_1, $$
where $S_\Ch = E \Ch^+ - \{ \bar{\gamma}_1,\dots,\bar{\gamma}_g\}$ for a collection of edges $\gamma_i$ outside a fixed spanning tree (so $S_\Ch$ consists of the edges of a spanning tree with both orientations, and one oriented edge for each geometrical edge outside that spanning tree).

As was observed by Jakub Byszewski, there is a natural injective group homomorphism
$ \varphi : Z_1 \hookrightarrow \ker(1-T)$ given as follows: if a cycle $[c]$ is represented as $\displaystyle \sum_{e \in I} e$, we define \begin{equation} \label{defphi} \varphi([c]) = \sum_{e \in I} e - \sum_{e \in I} \bar{e}. \end{equation}
This does not depend upon the choice of a representative for the cycle $c$. 
Let us check that, indeed, $(T-1)\varphi([c])=0$. Fix $e \in I$, and suppose $e' \in I$ is the (unique) edge in $I$ such that $t(e)=o(e')$. Observe the basic equation $T(e)-e' = T(\bar{e}')-\bar{e},$ that is illustrated in Figure \ref{fig3}. 

\begin{figure}[h]
   \centering
   \input{fig3.pstex_t} 
   \caption{A basic equation}
   \label{fig3}
\end{figure}

We can therefore calculate
\begin{eqnarray*} 
(T-1)\varphi([c]) &=& (T-1) (\sum e - \sum \bar{e}) \\ &=& \sum T(e) - \sum e  - \sum T(\bar{e}) + \sum \bar{e} \\ &=& \sum T(\bar{e}') - \sum \bar{e} + \sum e' - \sum T(\bar{e}) - \sum e + \sum \bar{e} \\ & = & 0. \end{eqnarray*} 
The map $\varphi$ is injective: choose a basis for the space of loops that doesn't contain any edges in sinks, then if $e \in I$, $\bar{e} \notin I$, so the independence of $e$ and $\bar{e}$ in the image implies injectiveness.

Below is the promissed geometrical interpretation of the kernel of $T$:
 \end{se}

\begin{lem} \label{jakub} For $g \geq 2$, the kernel of $1-T$ is isomorphic to the cycle space $Z_1$ via the isomorphism $\varphi$ defined in \textup{(\ref{defphi})}: $$K_1(\cO_{\Ch}) = \ker(1-T) = Z_1 \cong \Z^g. $$ \end{lem}

\begin{proof}
From \ref{oli}, we know that the kernel of $1-T$ has the same rank ($=g$) as $Z_1$. Alternatively, this follows independently from computations with the graph zeta function as in \ref{graphzeta} below.
Since $\varphi$ is also injective, the image is of the form $\bigoplus\limits_{i=1}^g a_i \Z$ for some $a_i \in \Z_{>0}$. We should prove that $a_i=1$ for all $i$. 

An \emph{end} of $\Ch$ is a connected contractible subgraph of $\Ch$ that shares exactly one vertex with its complement in $\Ch$. 

Given any edge not belonging to an end, choose a loop $\gamma$ in which $e$ occurs with multiplicity one, and in which $\bar{e}$ doesn't occur. Then $e$ occurs in $\varphi(\gamma)$ with multiplicity one. 

Now suppose $e$ belongs to an end. Since $\varphi$ is bijective after tensoring with $\Q$, for any $\sum a_e e \in \ker(1-T)$ we can find $a,b \in \Z$ and a loop $\gamma$ such that 
$ a(\sum a_e e) = b \varphi(\gamma) = b (\sum e - \sum \bar{e})$. Since $e$ is in an end, it is not in the support of this last sum, so it cannot occur in $\sum a_e e$ either. 

The conclusion is that edges in ends don't occur in $\ker(1-T)$, and any other edge occurs with multiplicity one in some element of $\IM(\varphi)$. This would not be the case if some $a_i>1$.  
\end{proof}

\begin{rem}
For graphs without ends, this result is also found in \cite{Rob3}, where it is applied to prove the following: let $\cT$ denote the universal covering of $\Ch$ and $\Gamma$ its fundamental group. The group of $\Gamma$-invariant integral valued measures on clopen sets of the boundary $\partial \cT$ is isomorphic to $Z_1$.
\end{rem}

\begin{rem}
The map $\varphi$ gives a natural isomorphism $Z_1 \rightarrow \ker(1-T)$. On the other hand, Cuntz has constructed a natural isomorphism $\ker(1-T) \rightarrow K_1(\cO_\Ch)$ (cf.\ \cite{Rordam}, p.\ 33 and \cite{Cuntz}, 3.1). It would be interesting to give a direct formula for the map $Z_1 \rightarrow K_1(\cO_\Ch)$ that relates graph homology and operator $K_1$. From the proof of 3.1 in Cuntz \cite{Cuntz}, it is not hard to obtain an explicit form of the map $Z_1 \rightarrow K_1(\bar{\cO}_\Ch)$, where $\bar{\cO}_\Ch$ denotes the stabilization of $\cO_\Ch$. 
\end{rem}

\begin{se}[Relation to the Ihara graph zeta function] \label{graphzeta}
We can (independently of the above) compute the rank of the operator $1-T$ acting on the vector space $\Q^{E\Ch^+}$. For that, we use the relation between the characteristic polynomial of $T$ and the Ihara zeta function of the graph. Theorem 4.3 from \cite{Bass} (applied to the trivial representation) implies that $ \det(1-uT) = (1-u^2)^{g-1} \det(\Delta(u)) $
for $\Delta$ the graph Laplace operator. Since $g \geq 2$, the universal covering tree of $\Ch$ is not a linear tree, and \emph{loc.\ cit.}, Theorem 5.10.b(i) says that $\Delta(u)(1-u)D^+(u)$ with $D^+(1) \neq 0$, so 
$$ \ord_{u=1} \det(1-Tu) = g. $$ Now $x \in \ker(1-T)$ if and only if $x$ belongs to the eigenspace of $T$ for the eigenvalue $+1$, and we have just seen that this space is $g$-dimensional. Therefore, the rank of $1-T$ is $2m-g$. From $Z_1 \hookrightarrow \ker(1-T)$, and this result, one finds back Lemma \ref{jakub}. 

We can also reverse the logic and use the above (independent) calculation of the Smith Normal Form and kernel of $1-T$ to deduce facts about the Bass-Hashimoto $T$-operator and the graph zeta function. In particular, the statements about the multiplicity of the eigenvalue $+1$ for $T$ follow directly. 
\end{se}

\begin{figure}[h]
   \centering
   \input{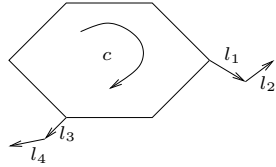} 
   \caption{Generators for $\ker(1-T)$ are $\varphi(c)$ and $c+\sum\limits_{i=1}^4 l_i$}
   \label{fige}
\end{figure}

\begin{rem}
If $g=1$, the universal covering of the graph is a linear tree with a cyclic action of $\Z$. This corresponds to Tate's uniformization of an elliptic curve with totally split multiplicative reduction. Then \cite{Bass}, 5.11.10 implies that $1-T$ has rank $2m-2$, which is compatible with \ref{oli}. It is easy to check the following: let $c$ denote a (fundamental) cycle in $\Ch$. Then two independent elements of $\ker(1-T)$ are given by $\varphi([c])$ and 
$ c  + \sum l$, where the sum is over all edges $l$ outside $c$ that point away from $c$, see Figure \ref{fige} for an example.  
\end{rem}

\begin{rem}
As noticed in the introduction, we can revert the logic of this section to apply operator algebra $K$-theory to graph theory as follows: given a  graph $\Ch$ of cyclomatic number $g \geq 2$, write $\Ch = \Gamma \backslash \cT$, where $\cT$ is the universal covering tree of $\Ch$ and $\Gamma$ is a free group on $g$ generators. Then Kumjian and Pask (\cite{KP}) have shown that $\cO_\Ch \cong C(\Lambda_\Gamma) \rtimes \Gamma$ in the sense of strong Morita equivalence. By the Brown-Green-Rieffel result from \ref{BGR}, we find that the stable isomorphism class of $\cO_\Ch$ only depends on $g$.

In \cite{Rob}, Robertson showed  that $C(\Lambda_\Gamma) \rtimes \Gamma$ is stricly isomorphic to the Cuntz-Krieger algebra associated to the stable graph in Figure \ref{fig1}. We can therefore compute $K_0(\cO_\Ch)$ for any chosen $\Ch$, and doing this for the ``flower'' implies that $K_0(\cO_\Ch)$ is as expected. This, in its turn, implies the results about the image and rank of $1-T$, independently of graph theoretical considerations.  
\end{rem}

\begin{rem} The proof in \cite{Rob} includes the calculation of the $K$-theory of the Cuntz-Krieger algebra for Figure \ref{fig1}, but directly ``on the boundary'', whereas our proof takes place on the reduction graph. This can be seen as another manifestation of a ``holography principle'', by which information on the boundary can be equivalently expressed on the tree itself, cf.\ \cite{MM}. It would be very interesting to extend this kind of calculation to higher dimensional buildings, cf.\ \cite{CMRV}, \cite{Rob2}. \end{rem}

\begin{se}[Analogue for Riemann surfaces] \label{co}
The analogue of this result in the classical theory of global uniformization of Riemann surfaces is as follows.  Let $A=C(\P^1(\R)) \rtimes \Gamma$, with $\Gamma$  a torsion-free cocompact lattice in $\PSL(2,\R)$, and $g$ the genus of the Riemann surface uniformized by $\Gamma$. Then Anan\-tha\-ra\-man-De\-la\-ro\-che (\cite{AD}) proved that the $K$-theory of $A$ is given by $K_0(A)=\Z^{2g+1} \oplus \Z/(2g-2)\Z$. The proof is topological (via a Thom isomorphism), and although it follows that $A$ is isomorphic to a Cuntz-Krieger algebra, there is no apparent direct link between the matrix of that algebra and some combinatorial structure on the Riemann surface, as is the case for the non-archimedean theory.
\end{se}

\begin{se}[Analogue for Mumford curves] \label{stable}
Let $k$ be a non-archimedean complete discretely valued field of mixed characteristic with absolute value $| \cdot |$. A projective curve $X$ over $k$ is called a \emph{Mumford curve} if
it is uniformized over $K$ by a \emph{Schottky group}. This means that there exists
a free subgroup $\Gamma$ in $PGL(2,K)$ of rank $g$, acting on ${\bf P}^1_K$ with \emph{limit set} $\Lambda_\Gamma$ such that $X$ satisfies $X^\an \cong \Gamma
\backslash ({\bf P}^{1,\an}_K - {\Lambda}_\Gamma)$ as rigid analytic
spaces. Mumford (\cite{Mumford}) has shown that these conditions
are equivalent to the existence of a stable model of $X$ over the ring of integers $\cO_k$ of $k$
whose special fiber consists only of rational components with double points over the residue field. 

Suppose the ground field $k$ is large enough so that the group $\Gamma$ acts on the Bruhat-Tits tree of $PGL(2,k)$ without inversions. This is always possible by a finite extension of $k$ if necessary. Let $\cT_\Gamma$ denote the subtree of the Bruhat-Tits tree spanned by geodesics connecting fixed points of hyperbolic elements in $\Gamma$. Then $\Ch_X:=\T_\Gamma/\Gamma$ is a finite graph that is intersection dual to the stable reduction of the curve $X$. In particular, the cyclomatic number of $\Ch_X$ equal the rank of $\Gamma$, which equals the genus $g$ of $X$ (cf.\ \cite{Mumford}, Thm.\ (3.3)). Note that $\Ch_X$ is allowed to have loops and multiple edges. 

As is observed in \cite{GvdP}, p.\ 124, any graph $\Ch$ can occur as the stable reduction graph of a Mumford curve, as soon as $\Ch$ is finite, connected, and every vertex which is not connected to itself is the origin of at least three edges. We call such a graph a \emph{stable graph}.

In \cite{ConsaniMarcolli}, a spectral triple was associated to $X$ in which the operator algebra is the boundary operator algebra of $\Gamma$. For this algebra, the Mumford curve plays the r\^ole of the Riemann surface in the results of \ref{co}. The current work is inspired by the question whether finer invariants of $X$ are detected by replacing this algebra with $\cO_\Ch$ for $\Ch=\Ch_X$ the stable reduction of $X$. 
\end{se}

\section{Strict isomorphism type}

\begin{se} The (not only stable) isomorphism type of $\cO_{\Ch}$ is determined by the image of the unit of that algebra in its $K_0$-group. Suppose that $\Ch$ and $\Ch'$ are two stable graphs with the same cyclomatic number $g \geq 2$. To ease notation, add a prime to any symbol pertaining to $\Ch'$. Since we have already shown that $\cO:=\cO_{\Ch}$ and $\cO':=\cO_{\Ch'}$ are stably isomorphic, to prove that they are actually isomorphic, by R{\o}rdam (\cite{Rordam}) it suffices to decide whether or not there exists an isomorphism of abelian groups $K_0(\cO) \cong K_0(\cO')$ that maps the class of the unit $1$ of $\cO$ to that of the unit $1'$ of $\cO'$. We can now recast this in terms of linear algebra as follows: in our setting, an isomorphism of $K_0$-groups is given by an automorphism of $\Z^g \oplus \Z/(g-1)\Z$. Now Cuntz's isomorphism $\varphi: K_0(\cO) \rightarrow \Z^n/(1-A^t)\Z^n$ is explicitly given by $$\varphi([1])=\varphi([\sum e_i e^*_i]) = \sum \varphi([e_i e^*_i])=(1,\dots,1).$$ 
Hence it suffices to check whether there is an isomorphism $\Z^n/(1-A^t) \cong \Z^{n'}/(1-(A')^t)$ that fixes the class of ${\1}=(1,\dots,1)$. 

We now look at two examples to show that the isomorphism type can indeed vary.

\end{se}

\begin{ex}
 For a flower as in Figure \ref{fig1}, any preimage of $\1$ by $1-T$ has $g-1$ in the denominator; this is easily seen, since the matrix $B$ from \ref{oli}.1 is invertible. Hence the image of $\1$ in $K_0(\cO_{\Ch})$ is an element of exact order $g-1$. 
\end{ex} 

\begin{ex}
For a graph consisting of two vertices that are connected by $g+1$ edges, the matrix of $1-T$ is a block matrix of the form $\tinymat -1_{g+1} B B -1_{g+1} $. It is easy to check that $\1$ is the image of $$( -1,\underbrace{\frac{2}{g-1},\dots,\frac{2}{g-1}}_g,\frac{g+1}{g-1},\underbrace{0,\dots,0}_g ) .$$
In particular, if $g$ is odd, $(g-1)/2 \cdot {\1} \in \IM(1-T)$, so the order of $\1$ in $K_0(\cO_{\Ch})$ divides $(g-1)/2$. 
One can check that if $g$ is odd, $\1$ has order exactly equal to $(g-1)/2$ in $K_0(\cO_{\Ch})$, whereas if $g$ is even, the order is exactly equal to $g-1$. 
\end{ex}

\begin{rem}
Connes showed in \cite{Connestrans}, Cor.\ 6.7, that in the setting of Riemann surfaces, and sticking to the notations of Remark \ref{co}, the class of the identity $\1$ in $K_0(A)$ has order $2g-2$, the Euler characteristic of the Riemann surface. This fits with the results of \cite{AD} cited in \ref{co}, in the sense that the class of $\1$ generates the subgroup $\Z/(2g-2)\Z$ of $K_0(A)$.
\end{rem}

\begin{se}[Proof of Theorem \ref{thm2}] Let $\lambda$ denote the minimal positive integer such that the equation $$ (1-A) \cdot x = \lambda \1 $$
has a solution in an integral vector $x \in \Z^{2m}$; then, if it exists, $\lambda$ is the exact order of $\1$ in $K_0(\cO_\Ch)$.  We know from \ref{oli}.2 that there exist $X,Y \in GL_{2m}(\Z)$ such that 
$$ X(1-A)Y = \diag(1,\dots,1,g-1,\underbrace{0,\dots,0}_g),$$
where $X$ (resp.\ $Y$) is given by the row (resp.\ column) operations performed in the course of proving \ref{oli}. In particular, if $(1-A) \cdot x = \lambda \1$, then \begin{equation} \label{eq1} \diag(1,\dots,1,g-1,\underbrace{0,\dots,0}_g) \cdot y = \lambda X \cdot \1 \end{equation} has an integral solution $y \in \Z^{2m}$. The   equations corresponding to the first $2m-g-1$ rows of equation (\ref{eq1}) obviously have a solution in integers for any integral $\lambda$. The only question that remains is whether the remaining $g+1$ equations have an integral solution.

We calculate these equations by finding the entries of $X \cdot \1$ by performing the row operations speficied by $X$ on $\1$. At the start, $\1$ has entry $1$ at every place. The general Smith Normal Form reduction process in \ref{oli} starts in \ref{oli}.2. The row operations that where performed in \ref{oli}.2 correspond to the contraction of all non-loop edges, so that in the end, only one vertex remains. This operation can be seen as collapsing all vertices to a given one, i.e., to reach the final form of the matrix, we have to collapse $|V\Ch|-1$ vertices. To each of the contractions of an edge corresponds adding that row to any ingoing edge of the source of that edge. Hence each of the edges (then loops) that remain after the complete contraction process in \ref{oli}.2 has $|V\Ch|-1$ rows added to it. The result of this operation on $\1$ is that it has been tranformed into $$ \1 \leadsto (\dots, \underbrace{|V\Ch|, \dots, |V\Ch|}_{2g}).$$ 

Then the row operations in reducing the matrix from \ref{oli}.1 correspond to subtracting from the last $g$ rows the corresponding of the first $g$ rows, and then adding all rows to the $(2m-g)$-th one. 
After the reduction in \ref{oli}.1, one finally arrives at 
$$ \1 \leadsto X \cdot \1 = (\dots, g \cdot |V\Ch|, \underbrace{0, \dots, 0}_g). $$  
The last $g$ equations in (\ref{eq1}) therefore also admit a solution.
Now the final equation to consider is that on the $(g+1)$-to-last row:
$$ ? \ \exists z \in \Z \ : \ (g-1) z = \lambda \cdot g \cdot |V\Ch|, $$
and one sees that the minimal integral $\lambda$ for which a solution exists is $$ \lambda = \frac{g-1}{(g-1,|V\Ch|)}.$$
This proves the first part of the theorem, and proves in particular that $g-1$ anihilates the image of $\1$, that hence has finite order.   

There is an automorphism of the group $\Z^g \oplus \Z/(g-1)\Z$ that carries an element of $\Z/(g-1)\Z$ to another exactly if these elements have the same order in $\Z/(g-1)\Z$.

To see that any isomorphism type occurs for a stable graph $\Ch$ of genus $g \geq 2$, we will show that one can realise any number of vertices $m$ with $1 \leq m \leq g-1$ for such a stable graph. We will show that one can find a stable graph $\Ch_g$ of genus $g$ with $2g-2$ vertices. Collapsing these vertices one after the other keeps the graph stable of the same genus, decreasing the number of vertices by one every time. The graph $\Ch_g$ is as in Figure \ref{figs}. It has $2g-2$ vertices. The two terminal vertices have a loop attached to it, and then consecutive vertices are connected alternatingly by a simple edge and by two edges. The genus of $\Ch_g$ is easily seen to be $g$, and $\Ch_g$ is stable. 

 \begin{figure}
   \centering
   \includegraphics{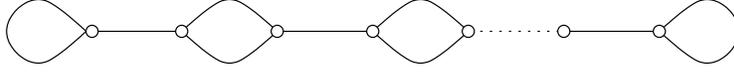} 
   \caption{The stable graph $\Ch_g$ with $2g-2$ vertices}
   \label{figs}
\end{figure}
Since for fixed $g$, all contractions of $\Ch_g$ are homotopic but the isomorphism type of the corresponding operator algebra varies, strict isomorphism is not a homotopy invariant.

Since the boundary operator algebra of genus $g$ is strictly isomorphic to $\cO_{\Ch}$ with $\Ch$ a ``flower'' as in Figure \ref{fig1} (see Robertson \cite{Rob}) and in $K_0$ of the latter algebra, $\1$ has order $g-1$, it follows that only algebras $\cO_\Ch$ for graphs with $$\frac{g-1}{(g-1,|V\Ch|)}=g-1,$$ i.e., graphs for which $(g-1,|V\Ch|)=1$, are strictly isomorphic to boundary operator algebras. This finishes the proof of Theorem \ref{thm2}. \hfill $\Box$
\end{se}

\begin{rem}
Note that by Euler's formula, $$ \mbox{ord } \1 =  \frac{g-1}{(g-1,|V\Ch|)} = \frac{g-1}{(g-1,|E\Ch|)}.$$
\end{rem}


\begin{small}

\end{small}


\finalinfo 

\end{document}